\documentclass[11pt]{amsart}
\usepackage{mathrsfs}
\usepackage{amssymb,latexsym}

\usepackage{amssymb,amsmath,ams}
\newdimen\AAdi%
\newbox\AAbo%
\def\AAk#1#2{\s_etbox\AAbo=\hbox{#2}\AAdi=\wd\AAbo\kern#1\AAdi{}}%
\def\AAr#1#2#3{\s_etbox\AAbo=\hbox{#2}\AAdi=\ht\AAbo\raise#1\AAdi\hbox{#3}}%
\font\tenmsb=msbm10 at 11pt \font\sevenmsb=msbm7 at 8pt
\font\fivemsb=msbm5 at 6pt
\newfam\msbfam
\textfont\msbfam=\tenmsb \scriptfont\msbfam=\sevenmsb
\scriptscriptfont\msbfam=\fivemsb

\newlength{\abstractwidth}
\newcommand{\be}{\begin{equation}}
\newcommand{\bea}{\begin{eqnarray}}
\newcommand{\eea}{\end{eqnarray}}
\newcommand{\ee}{\end{equation}}
\newcommand{\<}{\langle}
\renewcommand{\>}{\rangle}
\def\ba{\begin{eqnarray}}
\def\ea{\end{eqnarray}}

\def\14{{1 \over 4}}

\def\p{\partial}
\def\n{\nabla}

\def\e{\varepsilon}
\def\a{\alpha}
\def\b{\beta}

\def\g{\gamma}

\def\vp{\varphi}

\def\log{{\rm log}}

\newtheorem{theorem}{Theorem}

\newtheorem{proposition}{Proposition}
\newtheorem{lemma}{Lemma}

\begin{document}

\newcommand{\irrep}[1]{${\bf {#1}}$}

\baselineskip=15pt
\setcounter{equation}{0}
\setcounter{footnote}{0}

\title[A maximum rank problem]
{A maximum rank problem for degenerate elliptic fully nonlinear equations}
\author{Pengfei Guan}
\address{Department of Mathematics\\
         McGill University\\
         Montreal, Quebec. H3A 2K6, Canada.}
\email{guan@math.mcgill.ca}
\author{D. H. Phong}
\address{Department of Mathematics\\
Columbia University\\
New York, NY 10027, USA.}
\email{phong@math.columbia.edu }

\thanks{Research of the first author was supported in part
by NSERC Discovery Grant. Research of the second author was
supported in part by the National Science Foundation
grant DMS-07-57372.}

\begin{abstract}
The solutions to
the Dirichlet problem for two degenerate elliptic fully nonlinear equations
in $n+1$ dimensions, namely the real Monge-Amp\`ere equation and the Donaldson equation, are shown to have maximum
rank in the space variables when $n \leq 2$. A constant rank property is also
established for the Donaldson equation when $n=3$.
\end{abstract}

\maketitle

\section{Introduction}
\setcounter{equation}{0}

The solutions of elliptic partial differential equations are known to have many remarkable convexity
properties, under suitable structure conditions.
Some early works are those of Brascamp-Lieb \cite{BL}, Caffarelli and Friedman \cite{CF},
Yau \cite{Y}, with many important subsequent developments
(see \cite{ALL, KL, CGM, GLZ, L, BG1, BG2} and also references therein).
The constant rank theorem has been established for a general class of fully elliptic nonlinear equations. But the situation for degenerate elliptic fully nonlinear
equations has remained largely unexplored, despite its considerable interest
for example in geometry. One exception is the beautiful work of Lempert \cite{Lm} on the solution to the homogeneous
complex Monge-Amp\`ere equation on convex domains in ${\bf C}^n$ with prescribed $\log$ singularity at an interior point (the pluri-Green's function). Using a complex foliation, he showed that the solution is smooth and the complex Hessian has maximum rank $n-1$. Even for that result, there is no known PDE proof.

\smallskip
In this paper, we study a maximum rank problem for the Dirichlet problem for two basic models of such equations, on the space $X^n\times T$, where $X^n=({\bf R}/{\bf Z})^n$
is the $n$-dimensional torus and $T=(0,1)$ is the unit interval.
The first model is the Monge-Amp\`ere equation
\bea
\label{MAeq}
{\rm det}(D_{x,t}^2 u+I_{n+1})=\e
\eea
and the second is the equation
introduced by Donaldson \cite{D2}
\bea
\label{Deq}
u_{tt}(n+\Delta u)-\sum_{j=1}^n u_{jt}^2=\e.
\eea
Here the variables in $X\times T$
have been denoted by $(x,t)$, and $I_{n+1}$ is the $(n+1)\times (n+1)$ matrix
with the $n\times n$ identity matrix $I_n$ as its upper left block, and all zeroes
on its $(n+1)$ row and its $(n+1)$ column.
The right hand side $\e$ is a strictly positive constant, but may be arbitrarily
small. One imposes the
Dirichlet condition
\bea
\label{boundarycondition}
u(x,0)=u^0(x),
\qquad u(x,1)=u^1(x),
\eea
where the boundary data $u^0$ and $u^1$ are assumed to be $C^\infty$.
For the Monge-Amp\`ere equation, the solution $u$ is required to satisfy
$D_{xt}^2u+I_{n+1}\geq 0$, while for the Donaldson equation,
it is required to satisfy $n+\Delta u\geq 0$.

\smallskip

Both cases of right hand side $0<\e<<1$ and
$\e=0$ are of importance in geometry and physics:
the geodesic and approximate geodesic equations for the space of K\"ahler potentials
coincide with a complex version
of (\ref{MAeq}) \cite{D1,M,S}, and for toric varieties, they reduce
to a real version of (\ref{MAeq}) on polytopes with Guillemin boundary conditions.
The equation (\ref{Deq}) has a similar interpretation as the geodesic and approximate geodesic equations for the space of volume forms on a Riemannian
manifold. It coincides with (\ref{MAeq}) when $n=1$, but is more closely related
to Nahm's equation in theoretical physics as well as to some free boundary problems in applied mathematics \cite{D2}.

\smallskip
For fixed $\e>0$, the equations (\ref{MAeq})
and (\ref{Deq}) are elliptic. The existence of a unique smooth solution for $\e>0$
is a consequence of the general theory \cite{CNS, Gb}
for the Monge-Amp`ere equation (\ref{MAeq}),
and it has been established in \cite{CH, H} for the Donaldson equation
(\ref{MAeq}). For $\e=0$, it
has been shown by D. Guan \cite{Gd}
that the partial Legendre transform of the solution of (\ref{MAeq})
is a linear function of $t$, and thus the equation admits a smooth solution
which is strictly convex in $x$. The existence of $C^{1,\a}$ solutions for
(\ref{Deq}) is in \cite{CH,H}. We refer to
\cite{C, PS1, PS2, PS3} for various regularity results
on the complex
Monge-Amp\`ere equation on K\"ahler manifolds.

\smallskip
The main question of interest is whether the equations (\ref{MAeq})
and (\ref{Deq}) have maximum rank in the space directions, in the sense that the Hessian satisfies
$D_x^2u+I_n>0$ for all $t$. Of course, the boundary data have to satisfy
the maximum rank property first. That is, for some strictly positive constant
$\lambda$,
\bea
\label{strictconvexity}
D_x^2u^0+I_n\geq \lambda,
\qquad
D_x^2u^1+I_n\geq\lambda.
\eea
We shall be interested in when there is an estimate
\bea
\label{maximumrank}
D_x^2u+I_n\geq \widetilde\lambda,
\eea
with some constant $\widetilde\lambda>0$ {\it uniformly in $\e$}? It does not appear that the
Monge-Amp\`ere and Donaldson equations (\ref{MAeq}, \ref{Deq}) fall under the scope of the broad
structure conditions which have been introduced for partial convexity properties
in \cite{BG2}. However, by building on the strong maximum principle
methods of \cite{CF, Y, CGM, BG1, BG2},
exploiting the specific form of the Monge-Amp\`ere and the Donaldson
equations, and pushing the desired estimates to the
boundary, we can establish the following:

\begin{theorem}
\label{1}
Let $u$ be the solution of the Monge-Amp\`ere equation (\ref{MAeq}) on $X^n\times T$ with $D_{xt}^2u+I_{n+1}\geq 0$ and
Dirichlet data (\ref{boundarycondition}) satisfying
the strict
convexity condition (\ref{strictconvexity}). Assume that $n\leq 2$.
Then for all $t\in T$
and all $\e>0$, $u(\cdot,t)$ satisfies the same strict convexity condition
in the interior,
\bea
\label{maximumrank1}
D_x^2u(x,t)+I_n\geq\lambda,
\eea
with the same $\lambda>0$ as in (\ref{strictconvexity}).
\end{theorem}

When $\e=0$, the solutions of (\ref{MAeq}) are explicit and manifestly
satisfy the inequality
(\ref{maximumrank1}) \cite{Gd}, so the interest in Theorem \ref{1} lies in
the solutions for $\e>0$ themselves. They can be easier to use than
the solutions for $\e=0$, see for example the complex case treated in
\cite{C, CS, B}. For the equation (\ref{Deq}), lower bounds
for $D_x^2u+I_n$ in both
cases $\e>0$ and $\e=0$ were not known. We have

\begin{theorem}
\label{2}
Let $u$ be the solution of the Donaldson equation (\ref{Deq}) on $X^n\times T$ with $n+\Delta u\geq 0$ and Dirichlet data (\ref{boundarycondition}) satisfying
the strict
convexity condition (\ref{strictconvexity}). Assume that $n\leq 2$.
Then the strict convexity condition (\ref{strictconvexity})
with the same lower bound $\lambda$ is preserved in the interior,  that is, for all $t\in T$
and all $\e>0$,
$u(\cdot,t)$ satisfies
\bea
D_x^2u(x,t)+I_n\geq\lambda.
\eea
\end{theorem}

It would be interesting to determine whether one can lift the restriction of $n\le 2$ in Theorem \ref{1} and Theorem \ref{2}. For the Donaldson equation (\ref{Deq}), we can prove the following partial constant rank theorem.

\begin{theorem}
\label{th2}
Suppose $\Omega$ is a domain in ${\bf R}^n$ and $\delta>0$. Let $u$ be a solution of the Donaldson equation (\ref{Deq}) on $\Omega \times (0, \delta)$ satisfying $D_{x}^2u+I_n\geq 0$ for each $t\in (0,\delta)$. Assume that $n\leq 3$.
Then the rank of $(D_{x}^2u+I_n)$ is constant for all $(x,t)\in \Omega \times (0, \delta)$.
\end{theorem}

For the Monge-Amp\`ere equation, $\e>0$ and the convexity condition
$D_{xt}^2u+I_{n+1}\geq 0$ imply trivially the strict convexity condition
$D_{xt}^2u+I_{n+1}>0$ and hence the strict space convexity
condition $D_x^2u+I_n>0$. Thus the main interest in
Theorem \ref{1} lies in the fact that the lower bound for
$D_x^2u+I_n$ depends only on the boundary data.
For the Donaldson
equation, even the mere space convexity of the solution does
not seem so easy.
We observe that it
follows from Theorem \ref{th2} when $n\leq 3$:

\begin{theorem}
\label{th3}
Let $u$ be the solution of the Donaldson equation (\ref{Deq}) on $X^n\times T$ with $n+\Delta u\geq 0$ and Dirichlet data (\ref{boundarycondition}) satisfying
the strict
convexity condition (\ref{strictconvexity}). Assume that $n\leq 3$.
Then
the strict convexity of $u$ is preserved in the interior,
that is, for all $t\in T$ and all $\e>0$,
\bea
D_x^2u(x,t)+I_n>0.
\eea
\end{theorem}

\medskip

There are many important questions related to the maximum rank problem
which should be investigated. Perhaps of greatest interest
is the question of whether maximum rank theorems such
as Theorem \ref{1} hold for the complex Monge-Amp\`ere
equation, i.e., the geodesic and approximate geodesic equations for the space of
K\"ahler metrics. It is not clear whether the techniques developed in \cite{L,GLZ} for complex nonlinear equations can be adapted to treat the maximum rank problem for the complex Monge-Amp\`ere equation. One would also
like to generalize Theorem \ref{1} to general
Riemannian manifolds of arbitrary dimension. The results of this paper
should be thought of as experimental. It is our hope that the paper can generate some interest for the study of the maximum rank problem, as we believe that it is an important topic in PDE and differential geometry.

\medskip

The proof of Theorems \ref{1}, \ref{2} and \ref{th2} is given in \S 5.
The essential part is contained in Propositions \ref{MAresult} and
\ref{Dresult}, which are proved in \S 3 and \S 4 respectively.

\section{The general set-up}
\setcounter{equation}{0}

Both the Monge-Amp\`ere and the Donaldson equations (\ref{MAeq})(\ref{Deq})
are equations of the form
\bea
\label{F}
F(D_{xt}^2u+I_{n+1})=\e,
\eea
where $F(M)$ is a function of
the symmetric $n\times n$ matrix $M=(M_{\a\b})$.
Recall that $u$ is a function on $(x,t)\in X^n\times T$. It is convenient
to denote by Latin letters $i,j,\cdots$ the $n$ indices
for the ``space" variables $x=(x_j)$, and by Greek letters $\a,\b,\cdots$ the
$n+1$ indices for the ``space-time" variables $(x,t)$.
As usual, we denote by $F^{\a\b}$ and $F^{\a\b,\g\delta}$ the derivatives of $F$ with respect to $M_{\a\b}$,
\bea
F^{\a\b}={\p F\over \p M_{\a\b}},
\qquad
F^{\a\b,\gamma\delta}={\p^2 F\over \p M_{\a\b} \p M_{\g\delta}}.
\eea

Let $\mu_0$ be the minimum over $X^n\times T$ of the lowest eigenvalue of
$D_x^2u+I_n$,
\bea
\label{mu0}
\mu_0={\rm min}_{(x,t)\in X^n\times T}{\rm min}_{|\xi|^2=1}
{\<(D_x^2 u(x,t)+I_n)\xi,\xi\>}.
\eea
We would like to show that $\mu_0$ is attained at the boundary. For this,
it suffices to show that the set where the matrix $D_x^2u+I_n-\mu_0 I_n$
has a zero eigenvalue is open. In practice, it suffices to show that for each
$K$, the set where the matrix $D_x^2u+I_n-\mu_0 I_n$ has a zero eigenvalue
of multiplicity $K$ is open. Let $x_0$ be an interior point of $X^n\times T$
where $D_x^2+I_n-\mu_0 I_n$ has a zero eigenvalue of multiplicity
$K$. Set
\bea
\label{K}
\vp=\sum_{i_1<\cdots <i_{n-K+1}}\lambda_{i_1}\cdots \lambda_{i_{n-K+1}}
\equiv
\sigma_{n-K+1}(\lambda_1,\cdots,\lambda_n)
\eea
where $\lambda_i$ are the eigenvalues of $D_x^2u+I_n-\mu_0I_n$
\footnote{The function $\vp$ depends obviously on the choice
of order $K$. To lighten the notation for
$\vp$, we have not indicated this explicitly.}.
The strong maximum principle
reduces the desired statement to a key local, elliptic
inequality near $x_0$. The precise formulation we need is
the following:

\begin{proposition}
\label{strategy}
{\rm (a)} Let $x_0$ be an interior point where
$D_x^2u+I_n-\mu_0I_n$ has a zero eigenvalue
of some multiplicity $K$, and let $\vp$ be defined as in (\ref{K}).
If there is a constant $C$ so that
\bea
\label{keyestimate}
F^{\a\b}\vp_{\a\b}
\leq C(\vp+|\n\vp|)
\eea
for all points in a neighborhood of $x_0$, then $\vp$ vanishes
in a neighborhood of $x_0$.

{\rm (b)} If (\ref{keyestimate}) holds for an arbitrary
interior point $x_0$ where $D_x^2u+I_n-\mu_0I_n$ has a zero
of multiplicity $K$, and if  $\vp$ vanishes at some point,
then $\vp$ vanishes identically on $X\times T$.
In particular,
$\lambda$ is the largest lower bound for the boundary data $D_x^2u^0+I_n$
and $D_x^2u^1+I_n$,
and we have $\mu_0=\lambda$ and, for all $(x,t)\in X^n\times T$,
\bea
D_x^2u+I_n\geq\lambda I_n.
\eea
\end{proposition}

\bigskip
Thus we need to investigate estimates of the form (\ref{keyestimate}).
Let $x_0$ be an interior point as in Proposition \ref{strategy}, (a),
and let $x$ be an arbitrary point in a neighborhood of $x_0$.
Choose a parametrization $\lambda_1,\cdots\lambda_n$ of the eigenvalues of the matrix
$D_x^2 u+I_n-\mu_0I_n$ which is continuous in a neighborhood of $x_0$.
For each fixed $x$ in this neighborhood, we can choose a coordinate system
with $D_x^2u$ diagonal at $x$. Thus
at $x$, we have
\bea
u_{ij}+(1-\mu_0)\delta_{ij}=\lambda_i\delta_{ij}.
\eea
Define the matrix $v_{ij}$ by
\bea
v_{ij}=u_{ij}+(1-\mu_0)\delta_{ij}.
\eea
We divide the indices $i$, $1\leq i\leq n$,
into two sets of indices,
\bea
\{1,\cdots,n\}=G\cup B
\eea
with the ``good" set $G$ consisting of those indices $i$ for which
$\lambda_i(x_0)\not=0$, and the ``bad" set $B$ consisting of those indices $i$ for which $\lambda_i(x_0)=0$. Note that $\#G=n-K$ and $\#B=K$,
where $\#G,\#B$ denote the cardinalities of $G$ and $B$.
The starting point of our considerations is the following

\begin{lemma}
\label{starting}
Let $\#G$ be the number of good directions, and set
$\vp=\sigma_{\#G+1}(\lambda_1,\cdots,\lambda_n)$. Then we have

{\rm (a)} The function $\vp$ is of size
\bea
c_1\sum_{m\in B}v_{mm}\leq \vp
\leq c_2\sum_{m\in B}v_{mm}
\eea
for some strictly positive constants $c_1,c_2$.

{\rm (b)} The first derivatives of $\vp$ are given by
\bea
\label{firstderivatives}
\vp_\a
=(\prod_{g\in G}v_{gg})\sum_{m\in B}u_{mm\a}+O(\vp)
\eea

{\rm (c)} The second derivatives of $\vp$ are given by
\bea
\label{secondderivatives}
\vp_{\a\b}
=(\prod_{g\in G}v_{gg})(\sum_{m\in B}u_{\a\b mm}
-
2\sum_{m\in B}\sum_{g\in G}\frac{u_{mg\a}u_{mg\b}}{v_{gg}})
+
O(\vp+|\n\vp|).
\eea

{\rm (d)} The linearized operator $F^{\a\b}\vp_{\a\b}$ is given by
\bea
\label{linearized}
F^{\a\b}\vp_{\a\b}
&=&
-(\prod_{g\in G}v_{gg})(\sum_{m\in B}F^{\a\b,\g\delta}u_{\a\b m}u_{\g\delta m}
+
2\sum_{m\in B}\sum_{g\in G}\frac{F^{\a\b}}{v_{gg}}u_{mg\a}u_{mg\b})
\nonumber\\
&&
+
O(\vp+|\n \vp|)
\eea
\end{lemma}

\noindent
{\it Proof}: The function $\vp$ is a linear superposition
of terms, each of which is a product of $\#G+1$ eigenvalues
of $D_x^2u+I_n-\mu_0I_n$. Thus it is of the size of the sum
of the terms with exactly one eigenvalue in $B$. This establishes (a).
Formally, (b) and (c) can be established in the same way, by differentiation
of the eigenvalues if they are smooth. More generally, the same proof
can be adapted by expanding $\vp$ in terms of minors. As for (d),
successive differentiations of the equation $F(D_{xt}^2u+I_{n+1})=\e$ gives
\bea
&&
F^{\a\b}u_{\a\b \mu}=0\nonumber\\
&&
F^{\a\b}u_{\a\b \mu\mu}=-F^{\a\b,\g\delta}u_{\a\b\mu}u_{\g\delta\mu}.
\eea
Multiplying the expression for $\vp_{\a\b}$ in (c) by $F^{\a\b}$,
and making use of this last identity gives (d). Q.E.D.

\section{The Monge-Amp\`ere equation}
\setcounter{equation}{0}

We consider now more specifically the Monge-Amp\`ere equation, where
\bea
\label{MA}
F(M)={\rm det}\, M_{\a\b}
\eea
and $\e$ is a strictly positive constant.
Our main results in this situation can be stated as follows:

\begin{proposition}
\label{MAresult}
Let $u$ be a solution of the equation $F(D_{xt}^2u+I_n)=0$ on $X^n\times T$, which is convex in the sense that $D_{xt}^2u+I_{n+1}\geq 0$.
Define $\mu_0$ as in
(\ref{mu0}), and let $K$ be either $n$ or $n-1$.
Then the set of interior points $x_0$ where
the matrix $D_x^2u+I_n-\mu_0I_n$ has a zero eigenvalue
of multiplicity  $K$ is open.
\end{proposition}

\noindent
{\it Proof}: Let $x_0$ be an interior point where $D^2u+I_n-\mu_0I_n$ has a zero eigenvalue of multiplicity $K$.
We treat first the easier case when $K=n$.
In this case, the function $\vp$ is, explicitly,
\bea
\vp=\sum_{m=1}^n v_{mm}.
\eea
Consider next the expression (\ref{linearized}) for the function $F^{\a\b}\vp_{\a\b}$.
Since $|\n u_{ij}|\leq C\,\vp^{1\over 2}$ for $1\leq i,j\leq n$, we obtain modulo $O(\vp+|\n\vp|)$,
\bea
F^{\a\b}\vp_{\a\b}
=-\sum_{m\in B}F^{\a\b,\g\delta}u_{\a\b m}u_{\g\delta m}
=-2 \sum_{m\in B}F^{tt,\g\delta}u_{ttm}u_{\g\delta m}
=-2\sum_{m\in B}u_{ttm}\p_m F^{tt}.
\eea
For the Monge-Amp\`ere equation, $F^{tt}$ is simply the determinant
\bea
F^{tt}=\prod_{i=1}^n\mu_i.
\eea
where $\mu_i$ denotes the eigenvalues of $D^2u+I$.
Since we have then $\lambda_i=\mu_i-\mu_0$, we can write
\bea
F^{tt}
=\sum_{p=0}^n \sigma_{n-p}(\lambda_1,\cdots,\lambda_n)\mu_0^p.
\eea
Since $0\leq\lambda_i\leq\vp$ for all $i$, we have
\bea
\sum_{p=0}^{n-2}|\n \sigma_{n-p}(\lambda_1,\cdots,\lambda_n)|\leq C\,\vp,
\qquad
|\n\sigma_1(\lambda_1,\cdots,\lambda_n)|
=|\n\vp|.
\eea
Thus $|\p_mF^{tt}|\leq C(\vp+|\n\vp|)$.
Altogether, we obtain the inequality (\ref{keyestimate})
and the desired statement follows from Proposition \ref{strategy}.

\medskip

We consider now the case when the matrix $D_x^2u+I_n-\mu_0I_n$
admits at an interior point $x_0$ a zero eigenvalue
of multiplicity $n-1$. Thus there is only one good direction, which
we label $g$, $G=\{g\}$, and all other space directions,
$B=\{1\leq m\leq n;\  m\not=g\}$ are bad.
The expression (c) for $F^{\a\b}\vp_{\a\b}$ in Lemma \ref{starting}
becomes in this case
\bea
F^{\a\b}\vp_{\a\b}
=
-\sum_{m\in B}
\big(v_{gg}\sum_{\a\b,\g\delta}F^{\a\b,\g\delta}v_{\a\b m}v_{\g\delta m}
+2\sum_{\a\b}F^{\a\b}u_{mg\a}u_{mg\b})+
O(|\n\vp|+|\vp|).
\eea
The next step is
to derive an identity for the term $F^{\a\b}u_{mg\a}u_{mg\b}$:

\begin{lemma}
We have
\bea
F^{\a\b}u_{mg\a}u_{mg\b}
=
(u_{gg}+1)\sum_{\a,\b\not=g}F^{\a\b,gg}
(u_{mg\a}u_{mg\b}-u_{mgg}u_{\a\b m})
+
O(|\n\vp|+|\vp|)
\eea
\end{lemma}

\noindent
{\it Proof:} Recall that differentiating the equation gives
\bea
F^{\a\b}u_{\a\b m}=0
\eea
Extracting the terms involving the good direction $g$ gives,
\bea
F^{gg}u_{ggm}+2\sum_{\a\not=g}F^{\a g}u_{\a gm}
=
-
\sum_{\a,\b\not=g}F^{\a\b}u_{\a\b m}.
\eea
Returning to the expression $F^{\a\b}u_{mg\a}u_{mg\b}$, we can write
\bea
F^{\a\b}u_{mg\a}u_{mg\b}
&=&
F^{gg}u_{mgg}u_{mgg}+2\sum_{a\not=g}F^{\a g}u_{mg\a}u_{mgg}
+
\sum_{\a,\b\not=g}F^{\a\b}u_{mg\a}u_{mg\b}
\nonumber\\
&=&
u_{mgg}(F^{gg}u_{mgg}+2\sum_{\a\not=g}F^{\a g}u_{mg\a})
+
\sum_{\a,\b\not=g}F^{\a\b}u_{mg\a}u_{mg\b}
\nonumber\\
&=&
\sum_{\a,\b\not=g}F^{\a\b}(u_{mg\a}u_{mg\b}-u_{mgg}u_{\a\b m}).
\eea
We exploit the fact that $F$ is an affine function of any of the entries to write
\bea
F^{\a\b}=F^{\a\b,gg}(u_{gg}+1)+F^{\a\b}_{\vert_{u_{gg}+1=0}}.
\eea
The identity in the lemma follows then from the following claim
\bea
\sum_{m\in B}\sum_{\a,\b\not=g}
F^{\a\b}_{\vert_{u_{gg}+1=0}}u_{mg\a}u_{mg\b}
&=&O(|\n\vp|+\vp)
\nonumber\\
\sum_{m\in B}\sum_{\a,\b\not=g}
F^{\a\b}_{\vert_{u_{gg}+1=0}}
u_{mgg}u_{\a\b m}&=&
O(|\n\vp|+\vp).
\eea
To see the first identity above, we note that the terms with $\a\in B$ and $\b\in B$
are $O(\vp)$. Thus we need only consider the terms with at least $\a$ or $\b$ equal
to $t$. But then the cofactor $F^{\a\b}_{\vert_{u_{gg}+1=0}}$ has either a full column
or a full row of zeroes, and must be $0$.

Next, we consider the second identity above. For the same reason as above,
$F^{\a\b}_{\vert_{u_{gg}+1=0}}=0$ if either $\a$ or $\b$ is equal to $t$.
Thus we can restrict to $\a,\b\in B$.
Write now
\bea
\sum_{\a,\b\not=g}F^{\a\b}_{\vert_{u_{gg}+1=0}}u_{\a\b m}
&=&
\sum_{\a,\b\in B}F^{\a\b}_{\vert_{u_{gg}+1=0}}u_{\a\b m}
+O(|\n\vp|+\vp)
\nonumber\\
&=&
\sum_{i\in B}F^{ii}_{\vert_{u_{gg}+1=0}}u_{iim}
+
2\sum_{i,j\in B, i<j}
F^{ij}_{\vert_{u_{gg}+1=0}}u_{ijm}
+
O(|\n\vp|+\vp).
\eea
By inspection, we observe that

\smallskip
$\bullet$ If $i,j\in B$, then $F^{ij}_{\vert_{u_{gg}+1=0}}=0$ unless $i=j$.

\smallskip
$\bullet$ If $i\in B$, then $F^{ii}_{\vert_{u_{gg}+1=0}}=u_{gt}^2 \prod_{j\in B,j\not=i} (u_{jj}+1)$.

\smallskip
The last identity implies
\bea
\sum_{m\in B}
F^{ii}_{\vert_{u_{gg}+1=0}}u_{iim}
&=&
u_{gt}^2\mu_0^{n-2}\sum_{m\in B}u_{iim}+O(\vp)
\nonumber\\
&=& O(|\n\vp|+\vp).
\eea
The lemma is proved. Q.E.D.

\medskip
For our purposes, it is convenient to rewrite the identity in the preceding lemma
in the following form:
note that $u_{\a\b m}u_{\g\delta m}=O(|\nabla \vp|+\vp)$ if both $\a$ and $\b $ are in $B$.
Since neither of them is $g$, we can assume that at least one of them is $t$.
Thus

\begin{lemma}
\label{lemma1}
We have
\bea
F^{\a\b}u_{mg\a}u_{mg\b}
&=&
(u_{gg}+1)\,\big(F^{tt,gg}u_{tgm}^2
+
2\sum_{i\in B}F^{it,gg}
u_{mgi}u_{mgt}\big)
\nonumber\\
&&-
(u_{gg}+1)
\sum_{\a,\b\not=g}
F^{\a\b,gg}u_{mgg}u_{\a\b m}
+
O(|\n\vp|+|\vp|).
\eea
\end{lemma}

Our next task is to simplify the expression
\bea
\sum_{\a\b,\g\delta}F^{\a\b,\g\delta}u_{\a\b m}u_{\g\delta m}.
\eea
First, we isolate the contribution of the index $(gg)$,
which will cancel out with the corresponding term from the first
identity,
\bea
\sum_{\a\b,\g\delta}F^{\a\b,\g\delta}u_{\a\b m}u_{\g\delta m}
=
2 \sum_{\a\b}F^{\a\b,gg}u_{\a\b m}u_{ggm}
+
\sum_{(\a\b),(\g\delta)\not=(gg)}
F^{\a\b,\g\delta}u_{\a\b m}u_{\g\delta m}.
\eea
Next, we work out the remaining
contributions. For this, it is convenient to introduce the following sets of indices,
$\mathcal{A}=\{(gt),(tg),(tt)\}$ , and $\mathcal{B}={\mathcal{A}}^c$, so that
$(\a\b)$ is in $\mathcal{B}$ if and only if at least one of the indices $\a$ or $\b$
is in $B$.

\medskip
$\bullet$ If both $(\a\b)\in {\mathcal B}$ and $(\g\delta)\in{\mathcal B}$,
then $|u_{\a\b m}|+|u_{\g\delta m}|=O(\vp^{1\over 2})$, and thus
these contributions are $O(\vp)$ and can be neglected.

\medskip
$\bullet$ The contributions when both $(\a\b)$ and $(\g\delta )$ are in ${\mathcal A}$
can be worked out explicitly,
\bea
\sum_{(\a\b)\in{\mathcal A},(\g\delta)\in{\mathcal A}}
F^{\a\b,\g\delta}u_{\a\b m}u_{\g\delta m}
=
F^{tg,gt}u_{tgm}^2+F^{gt,tg}u_{tgm}^2=
-2 u_{tgm}^2 \prod_{j\in B}u_{jj}=-2F^{tt,gg}u_{tmg}^2.
\eea

\medskip
$\bullet$ The remaining contributions are
\bea
2\sum_{(\a\b)\in{\mathcal B},
(\g\delta)\in{\mathcal A}}
F^{\a\b,\g\delta}u_{\a\b m}u_{\g\delta m}.
\eea

To identify these terms, we divide the set ${\mathcal B}$ of indices $(\a\b)$ with at least one
index in $B$ into three mutually disjoint sets:

\medskip
${\mathcal B}_0=\{(\a\b); \a\in B,\quad \b\in B\}$

${\mathcal B}_1=\{(\a\b); \a\in\{g,t\},
\quad \b\in B\}$

${\mathcal B}_2=\{(\a\b); \a\in B,\quad \b\in \{g,t\}\}$

\medskip
The sum breaks up correspondingly
\bea
2\sum_{(\a\b)\in{\mathcal B},
(\g\delta)\in{\mathcal A}}
F^{\a\b,\g\delta}u_{\a\b m}u_{\g\delta m}
&=&
2\sum_{a=0,1,2}
\sum_{(\a\b)\in {\mathcal B}_a}
\big(F^{\a\b,tt}u_{\a\b m}u_{ttm}
\nonumber\\
&&
\qquad
\quad
+
(F^{\a\b,gt}+F^{\a\b,tg})u_{\a\b m}u_{gtm}\big).
\eea

\medskip
Each of these terms can now be worked out explicitly.
First, we have
\bea
\sum_{(\a\b)\in {\mathcal B}_1\cup {\mathcal B}_2}F^{\a\b,tt}u_{\a\b m}u_{ttm}=0
\eea
because we can see by inspection that $F^{\a\b,tt}$
is given then by a matrix with a column or a row of $0$ and hence must be $0$.

Next, we have
\bea
\sum_{(ij)\in {\mathcal B}_0}
F^{ij,tt}u_{ij m}u_{ttm}
=
\sum_{i\in B}u_{iim}u_{ttm}\,\prod_{j\not=i, j\in B} (u_{jj}+1).
\eea
This is because $F^{ij,tt}=0$ unless $i=j$, and the entries $F^{ii,tt}$
can be easily computed, giving the formula above. Since we can replace $u_{jj}+1$ by $\mu_0$
modulo $O(\vp)$, we obtain
\bea
\sum_{(ij)\in {\mathcal B}_0}
F^{ij,tt}u_{ijm}u_{ttm}
=
\mu_0^{n-2}\sum_{i\in B} u_{iim}u_{ttm}+O(\vp)=
O(|\n\vp|+\vp).
\eea

It remains only to determine the sum
\bea
\sum_{a=0,1,2}\sum_{(\a\b)\in {\mathcal B}_a}
(F^{\a\b,tg}+F^{\a\b,gt})u_{\a\b m}u_{gtm}.
\eea

Consider first the contributions from $(\a\b)\in {\mathcal B}_2$,
i.e. $(\a\b)=(\a g)$ or $(\a\b)=(\a t)$. Then it is clear that we obtain
\bea
(\a\b)=(\a g):&&
F^{\a g,tg}+F^{\a g,gt}=F^{\a g,gt} =-F^{\a t,gg}\nonumber\\
(\a\b)=(\a t):&&
F^{\a t,tg}+F^{\a t,gt}=F^{\a t,tg}=-F^{\a g,tt}=0.
\eea
Thus we find
\bea
\sum_{(\a\b)\in {\mathcal B}_2}
(F^{\a\b,tg}+F^{\a\b,gt})u_{\a\b m}u_{gtm}
=-\sum_{\a\in B}F^{\a t,gg}u_{\a gm}u_{gtm}.
\eea
Similarly, the contributions from $(\a\b)\in{\mathcal B}_1$ correspond to $(\a\b)=(g\b)$ or $(\a\b)=(t\b)$, and work
out to be
\bea
(\a\b)=(g\b):&&
F^{g\b,tg}+F^{g\b,gt}
=F^{g\b,tg}=-F^{gg,t\b}
\nonumber\\
(\a\b)=(t\b):&&
F^{t\b,tg}+F^{t\b,gt}=-F^{tt,g\b}=0,
\eea
and
\bea
\sum_{(\a\b)\in{\mathcal B}_1}
(F^{\a\b,tg}+F^{\a\b,gt})u_{\a\b m}u_{gtm}=
-
\sum_{\b\in B}F^{gg,t\b}u_{g\b m}u_{tgm}.
\eea
Finally, we come to the contributions from $(\a\b)\in {\mathcal B}_0$.
Here it is seen by inspection that only $(\a\b)=(\a\b)$ will contribute,
and thus
\bea
\sum_{(\a\b)\in{\mathcal B}_0}
(F^{\a\b,tg}+F^{\a\b,gt})u_{\a\b m}u_{gtm}
&=&
\sum_{\a\in B}
(F^{\a\a,tg}+F^{\a\a,gt})u_{\a\a m}u_{gtm}
\nonumber\\
&=&
2\sum_{\a\in B}F^{\a\a,tg}u_{\a\a m}u_{gtm}.
\eea
But an inspection shows that for $\a=i\in B$,
\bea
F^{ii,tg}=u_{gt}\,\prod_{j\not=i,j\in B}(u_{jj}+1)
\eea
so that
\bea
\sum_{(\a\b)\in{\mathcal B}_0}
(F^{\a\b,tg}+F^{\a\b,gt})u_{\a\b m}u_{gtm}
&=&
2u_{gt}\sum_{i\in B}u_{iim}u_{gtm}\,\prod_{j\not=i,j\in B}(u_{jj}+1)
\nonumber\\
&=&
2u_{gt}\mu_0^{n-2}\sum_i u_{iim}u_{gtm}+O(\vp)\nonumber\\
&=& O(|\n\vp|+\vp).
\nonumber\\
\eea
In summary, we have proved

\begin{lemma}
\label{lemma2}
We have the following identity
\bea
\sum_{\a\b,\g\delta }
F^{\a\b,\g\delta }u_{\a\b m}u_{\g\delta m}
&=&
2\sum_{\a\b}F^{\a\b,gg}u_{\a\b m}u_{ggm}
\nonumber\\
&&-2 (F^{tt,gg}u_{tmg}^2
+2\sum_{i\in B}F^{it,gg}u_{igm}u_{gtm})
+O(|\n\vp|+\vp).
\eea
\end{lemma}

Comparing the identities in Lemmas \ref{lemma1} and \ref{lemma2},
we obtain the main lemma in this case,

\begin{lemma}
\label{lemma3}
We have
\bea
F^{\a\b}u_{mg\a}u_{mg\b}
=-(u_{gg}+1)\,F^{\a\b,\g\delta }u_{\a\b m}u_{\g\delta m}
+
O(|\n\vp|+\vp).
\eea
\end{lemma}

\noindent
and, finally, since $v_{gg}=u_{gg}+1-\mu_0$,
\bea
F^{\a\b}\vp_{\a\b}
=
2\mu_0 \,F^{\a\b,\g\delta }u_{\a\b m}u_{\g\delta m}
+
O(|\n\vp|+\vp).
\eea
The first term on the right hand side is negative, modulo $O(\vp+|\n\vp|)$:
indeed, Lemma \ref{lemma3} shows that it is given by $-(u_{gg}+1)^{-1}
F^{\a\b}u_{mg\a}u_{mg\b}$. But $u_{gg}+1>0$ and,
in the case of the Monge-Amp\`ere equation, the matrix
$F^{\a\b}$ is just the matrix of minors of $D_{xt}^2u+I_{n+1}$,
which is positive.
Thus we obtain again the key estimate (\ref{keyestimate}). Q.E.D.

\section{The Donaldson equation}
\setcounter{equation}{0}

In our notation, the Donaldson equation (\ref{Deq})
is an equation of the form (\ref{F}), with $F(M)$ given by
\bea
\label{D1}
F(M)=M_{tt}(1+\sum_{j=1}^n M_{jj})-\sum_{j=1}^nM_{jt}^2.
\eea
We again consider the Dirichlet problem on the space $X^n\times T$, with
the usual boundary condition (\ref{boundarycondition}). Our main result is the following:

\begin{proposition}
\label{Dresult}
Let $u$ be a solution of the Donaldson equation $F(D^2u+I')=0$ on $X^n\times T$ satisfying $D_{x}^2u+I_n\geq 0$,
with $F$ as in (\ref{D1}).
Define $\mu_0\ge 0$ as in
(\ref{mu0}), and let $K$ be either $n$ or $n-1$.
Then the set of interior points $x_0$ where
the matrix $D_x^2u+I_n-\mu_0I_n$ has a zero eigenvalue
of multiplicity  $K$ is open. If $\mu_0=0$ in (\ref{mu0}), then the set of interior points $x_0$ where
the matrix $D_x^2u+I_n$ has a zero eigenvalue
of multiplicity  $K\ge n-2$ is open.
\end{proposition}

\noindent
{\it Proof}: As in the previous sections, we work at an arbitrary point $x$
in a neighborhood of a given point $x_0$ where the matrix $D^2u+I_n-\mu_0I_n$
admits a zero eigenvalue of multiplicity $K$. The three values $K=n, n-1,n-2$
correspond respectively to $\#G=0,1,2$, where $\#G$ is the number of good directions. The most difficult case is $\#G=2$ (corresponding to the case $K=n-2$). Thus we
write down the calculations for general $\#G$, and then specialize
to the cases of interest.

\medskip

We use the same notations as in Section \S 2.
If we use $\varphi=\sigma_{n-K+1}$ as before, and apply Lemma \ref{starting}
with the function $F(M)$ corresponding to Donaldson's equation,
we would find
\begin{eqnarray*}
F^{\a\b}\varphi_{\a\b}=-2\sum_{m\in B}\tilde Q_m
\,(\prod_{g\in G}v_{gg})+O(\varphi+|\nabla\varphi|),
\end{eqnarray*}
with
\begin{eqnarray*}
\label{d6-0}
\tilde Q_m=u_{ttm}\Delta u_m -\sum_{j\in G} u_{tjm}^2+u_{tt}\sum_{j,k\in G}\frac{u_{jkm}^2}{v_{jj}}+(n+\Delta u) \sum_{j\in G}\frac{u_{tjm}^2}{v_{jj}}-
2\sum_{j\in G}\sum_{k=1}^n\frac{u_{tk}u_{tjm}u_{jkm}}{v_{jj}}.
\end{eqnarray*}
We notice that there are linear terms of the form $\nabla u_{km}, k,m \in B$.
When $n\ge 3$ and $K=n-2$, these
linear terms cannot be bounded by $\varphi+|\nabla\varphi|$. To overcome this obstacle, we use instead
\begin{eqnarray}
\label{newvarphi}
\varphi=\sigma_{n-K+1}+q, \quad \mbox{where $q=\frac{\sigma_{n-K+2}}{\sigma_{n-K+1}}$.}
\end{eqnarray}
The regularity and strong concavity of $q$ was
proved in \cite{BG1}. Following the arguments in \cite{BG1, BG2} (e.g, see eq. (60) in \cite{BG2}, in our case, $F$ is independent of $\nabla u, u, x$) for $\varphi$ defined in (\ref{newvarphi}), we obtain
\begin{eqnarray}
\label{d5-0}
\sum_{\alpha,\beta=1}^N F^{\alpha\beta}\varphi_{\alpha\beta}
&= & -2\sum_{m\in B}[\sigma_l(G)+\frac{\sigma_{1}^{2}(B|m)-
\sigma_{2}(B|m)}{\sigma_{1}^{2}(B)}]\tilde Q_m + O(\varphi +\sum_{i,j\in B}|\nabla
u_{ij}|)\nonumber \\
&& -\sum_{\alpha,\beta=1}^N F^{\alpha\beta}[\frac{\sum_{i\in
B}V_{i\alpha}V_{i\beta}}{\sigma^3_{1}(B)} +
\frac{\sum_{i,j\in B,i\neq j}
u_{ij\alpha}u_{ji\beta}}{\sigma_{1}(B)}],
\end{eqnarray}
where
\begin{eqnarray}
\label{d6}
\tilde Q_m&=&u_{ttm}\Delta u_m -\sum_{j\in G} u_{tjm}^2+u_{tt}\sum_{j,k\in G}\frac{u_{jkm}^2}{v_{jj}}+(n+\Delta u) \sum_{j\in G}\frac{u_{tjm}^2}{v_{jj}}-
2\sum_{j\in G}\sum_{k=1}^n\frac{u_{tk}u_{tjm}u_{jkm}}{v_{jj}}\nonumber\\
&=& Q_m +Q_m^*,
\end{eqnarray}
\bea
\label{Q-D1a}
Q_m&=&u_{ttm}\Delta u_m-\sum_{k\in G}u_{ktm}^2
+u_{tt}\sum_{j,k\in G}{u_{mjk}^2\over {1+u_{jj}}}\nonumber\\
&& +
(n+\Delta u)\sum_{j\in G}{u_{mjt}^2\over {1+u_{jj}}}
-
2\sum_{j,k\in G}{u_{tk}u_{tjm}u_{mjk}\over {1+u_{jj}}},
\eea
\begin{eqnarray}\label{QQm}
Q_m^*=\sum_{m\in B}\sum_{j\in G}F^{\a\b}(\frac1{v_{jj}}-\frac1{1+u_{jj}})u_{mj\a}u_{mj\b},
\end{eqnarray}
and $\sigma_l(B)=\sum_{i\in B}v_{ii}$, $\sigma_{l}(B|m)=\sum_{i\neq m, i\in B}v_{ii}$,
\begin{equation}\label{def-V}
 V_{i\alpha}=u_{ii\alpha}\sigma_{1}(B)-u_{ii}\Big(\displaystyle\sum_{j\in
 B}u_{jj\alpha}\Big).
\end{equation}
The term $\sum_{i,j\in B}|\nabla
u_{ij}|$ in (\ref{d5-0}) can be controlled by $\varphi, |\nabla \varphi|$ and the last term in (\ref{d5-0})
in the same way as in \cite{BG1, BG2}. We obtain, for some $C>0$,
\begin{eqnarray}\label{d5}
\sum_{\alpha,\beta=1}^N F^{\alpha\beta}\phi_{\alpha\beta}
\le -2C\sum_{m\in B}(Q_m +Q_m^*) + O(\varphi +|\nabla
\varphi|).
\end{eqnarray}
Since $\mu_0\ge 0$,
we have $\frac1{v_{jj}}-\frac1{1+u_{jj}}\ge 0$. It is easy to see that
$Q_m^*$ in (\ref{QQm}) is nonnegative by the positivity of $(F^{\a\b})$.
Thus we would be done if we can show that $Q_m\geq 0$, modulo
$O(\vp+|\n\vp|)$.
Since the contributions of each index $m\in B$ are
entirely similar,  we can consider them individually.
To simplify the notation, we set $m=1$, and drop the subindex $m$ from $Q_m$.

\subsection{Using the equation}

Differentiating the equation gives
\bea
u_{tt1}=-u_{tt}{\Delta u_1\over n+\Delta u}+{2\over n+\Delta u}
\sum u_{tj}u_{tj1}.
\eea
Thus $Q$ can be written as
\bea
\label{Q-D1b}
Q= A+B+C
\eea
with $A,B,C$ defined by
\bea
A&\equiv&
-u_{tt}{(\Delta u_1)^2\over n+\Delta u}+u_{tt}\sum_{j,k\in G}{u_{1jk}^2\over 1+u_{jj}}
\nonumber\\
B&\equiv&
{2\over n+\Delta u}
\Delta u_1\,\sum u_{tj}u_{tj1}
-
2\sum_{j,k\in G}{u_{tk}u_{tj1}u_{jk1}\over 1+u_{jj}}
\nonumber\\
C&\equiv&
\sum_{j\in G}u_{1jt}^2({n+\Delta u\over 1+u_{jj}}-1).
\eea

\subsection{The $A$ terms}

The $A$ terms can be re-written as follows, modulo $O(\vp+|\n\vp|)$,
\bea
A=
{u_{tt}\over 2(n+\Delta u)}
\sum_{j,k\in G, j\not=k}
{{((1+u_{kk})u_{jj1}-(1+u_{jj})u_{kk1})^2}\over {(1+u_{jj})(1+u_{kk})}}
+
u_{tt}\sum_{j,k\in G,j\not=k}{u_{1jk}^2\over 1+u_{jj}}.
\eea
To see this, just write
\bea
-{(\Delta u_1)^2\over n+\Delta u}
+
\sum_{j,k\in G}{u_{1jk}^2\over 1+u_{jj}}
&=&
\sum u_{1jj}^2 ({1\over 1+u_{jj}}-{1\over n+\Delta u})
-
\sum_{j\not=k}{u_{jj1}u_{kk1}\over n+\Delta u}+
\sum_{j\not=k}
{u_{1jk}^2\over 1+u_{jj}}
\nonumber\\
&=&
\sum {u_{1jj}^2\over (1+u_{jj})(n+\Delta u)}\sum_{k\not=j}(1+u_{kk})
-
\sum_{j\not=k}{u_{jj1}u_{kk1}\over n+\Delta u}+
\sum_{j\not=k}{u_{1jk}^2\over 1+u_{jj}}
\nonumber\\
&=&
{1\over 2}\sum_{j\not=k}
{1\over n+\Delta u}
({1+u_{kk}\over 1+u_{jj}}u_{jj1}^2+{1+u_{jj}\over 1+u_{kk}}u_{kk1}^2)\nonumber\\
&&
-
{1\over n+\Delta u}\sum_{j\not=k}u_{jj1}u_{kk1}
+
\sum_{j\not=k}{u_{1jk}^2\over 1+u_{jj}}
\nonumber\\
&=&
{1\over 2(n+\Delta u)}
\sum_{j\not=k}
{{((1+u_{kk})u_{jj1}-(1+u_{jj})u_{kk1})^2}\over (1+u_{jj})(1+u_{kk})}
+
\sum_{j\not=k}{u_{1jk}^2\over 1+u_{jj}}.
\nonumber
\eea

\subsection{The $B$ terms}

The $B$ terms can be re-written as
\bea
B= -
\sum_{j\not=k}{u_{tj}u_{tj1}\over (1+u_{jj})(n+\Delta u)}
(v_{jj1}(u_{kk}+1)-v_{kk1}(u_{jj}+1))
-
2\sum_{j\not=k}u_{tk}u_{tj1}u_{1jk}{1\over 1+u_{jj}}.
\eea

To see this, we decompose the $B$ terms as follows
\bea
{2\over n+\Delta u}
\Delta u_1\sum_{k\in G}u_{tk}u_{tk1}=
{2\over n+\Delta u}
\sum u_{jj1}\sum_{k\in G}u_{tk}u_{tk1}=
B_1+B_2
\eea
and
\bea
-2\sum_{j,k\in G}
{u_{tk}u_{tj1}u_{jk1}\over 1+u_{jj}}=
B_3-2\sum_{j,k\in G,j\not=k}{u_{tk}u_{tj1}u_{jk1}\over u_{jj}}
\eea
with the terms $B_1,B_2,B_3$ defined by
\bea
B_1&=& 2\sum_k {u_{tk}u_{tk1}u_{kk1}\over n+\Delta u}
\nonumber\\
B_2&=& 2\sum_{j\not=k}{u_{tk}u_{tk1}u_{jj1}\over n+\Delta u}
\nonumber\\
B_3&=&-2\sum {u_{tk}u_{tk1}u_{kk1}\over 1+u_{kk}}.
\eea
The terms $B_1$ and $B_3$ can be regrouped as
\bea
B_1+B_3
&=&
-2\sum u_{tk}u_{tk1}u_{kk1}({1\over 1+u_{kk}}-{1\over n+\Delta u})
\nonumber\\
&=&
-2\sum u_{tk}u_{tk1}u_{kk1}{1\over (1+u_{kk})(n+\Delta u)}
\sum_{j\not=k}(1+u_{jj}),
\eea
and, combined with $B_2$, as
\bea
B_1+B_2+B_3
&=&
-2\sum u_{tk}u_{tk1}u_{kk1}{1\over (1+u_{kk})(n+\Delta u)}
\sum_{j\not=k}(1+u_{jj})
+
2\sum_{j\not=k}{u_{tk}u_{tk1}u_{jj1}\over n+\Delta u}
\nonumber\\
&=&
2
\sum_{j\not=k}{u_{tk}u_{tk1}\over (1+u_{kk})(n+\Delta u)}
(u_{jj1}(1+u_{kk})-u_{kk1}(1+u_{jj})).
\nonumber
\eea

\subsection{A second formula for $Q$}

For convenience, we write here the formula for $Q$ obtained in this manner
\bea
\label{Q-D1c}
Q&=& u_{tt}\big\{
{1\over 2(n+\Delta u)}
\sum_{j\not=k} {((1+u_{kk})u_{jj1}-(1+u_{jj})u_{kk1})^2\over (1+u_{jj})(1+u_{kk})}
+
\sum_{j\not=k}{u_{1jk}^2\over 1+u_{jj}}\big\}
\nonumber\\
&&
-
{2\,u_{tj}u_{tj1}\over (1+u_{jj})(n+\Delta u)}
(u_{jj1}(1+u_{kk})-(1+u_{jj})u_{kk1})
-
2\sum_{j\not=k}{u_{tk}u_{tj1}u_{1jk}\over 1+u_{jj}}
\nonumber\\
&&
+
\sum_{j\in G}u_{1jt}^2({n+\Delta u\over 1+u_{jj}}-1).
\eea

We can complete the square in $u_{kk}u_{jj1}-u_{jj}u_{kk1}$, and obtain
\bea
\label{Q-D1d}
Q&=&
{1\over 2}\sum_{j\not=k} \big\{({u_{tt}\over (1+u_{jj})(1+u_{kk})(n+\Delta u)})^{1\over 2}
((1+u_{kk})u_{jj1}-(1+u_{jj})u_{kk1}) \nonumber\\
&&
+
({u_{tk}u_{tk1}\over 1+u_{kk}}-{u_{tj}u_{tj1}\over 1+u_{jj}})
({(1+u_{jj})(1+u_{kk})\over u_{tt}(n+\Delta u)})^{1\over 2}
\big\}^2
\nonumber\\
&&
+D+ u_{tt}\sum_{j\not=k} {u_{1jk}^2\over 1+u_{jj}}
-2 \sum_{j\not=k}{u_{tk}u_{tj1}u_{jk1}\over 1+u_{jj}}
\eea
where we have introduced the $D$ terms
\bea
\label{D}
D&\equiv&
-
{1\over 2}\sum_{j\not=k}{(1+u_{jj})(1+u_{kk})\over u_{tt}(n+\Delta u)}
({u_{tk}u_{tk1}\over 1+u_{kk}}-{u_{tj}u_{tj1}\over 1+u_{jj}})^2
+C
\nonumber\\
&=&
-
{1\over 2}\sum_{j\not=k}{(1+u_{jj})(1+u_{kk})\over u_{tt}(n+\Delta u)}
({u_{tk}u_{tk1}\over 1+u_{kk}}-{u_{tj}u_{tj1}\over 1+u_{jj}})^2
+\sum_{j\in G} u_{1jt}^2 ({n+\Delta u\over 1+u_{jj}}-1).
\eea
When $\#G\leq 1$, the term $Q$ reduces to the term $D$,
and $D$ reduces in turn to a manifestly positive expression.
Thus the cases $\#G\leq 1$ are now proved.

\subsection{The $D$ terms}

In the rest of the proof, we will assume $\mu_0=0$ and $\#G=2$.

Expanding the squares in the $D$ terms gives
\bea
D&=&
-
{1\over 2}
\sum_{j\not=k} {1\over u_{tt}(n+\Delta u)}({1+u_{jj}\over 1+u_{kk}}u_{tk}^2 u_{tk1}^2
+{1+u_{kk}\over 1+u_{jj}}u_{tj}^2u_{tj1}^2)
\nonumber\\
&&
+
\sum_{j\not=k}{1\over u_{tt}(n+\Delta u)}u_{tk}u_{tk1}u_{tj}u_{tj1}
+\sum u_{1jt}^2{1\over 1+u_{jj}}\sum_{k\not=j}(1+u_{kk})
\nonumber\\
&=&
\sum u_{1jt}^2({1\over (1+u_{jj})}
\sum_{k\not=j}(1+u_{kk})
-
\sum_{k\not=j} {u_{tj}^2(1+u_{kk})\over u_{tt}(n+\Delta u)})
+
\sum_{j\not=k}{u_{tk}u_{tk1}u_{tj}u_{tj1}\over u_{tt}(n+\Delta u)}
\nonumber\\
&=&
\sum_{j\in G}u_{1jt}^2\sum_{k\not=j}{1+u_{kk}\over 1+u_{jj}}{u_{tt}(1+\Delta u)-u_{tj}^2\over u_{tt}(n+\Delta u)}
+
\sum_{j\not=k}{ u_{tk}u_{tk1}u_{tj}u_{tj1}\over u_{tt}(n+\Delta u)}.
\nonumber
\eea
We can now make use of the equation
\bea
u_{tt}(n+\Delta u)-u_{tj}^2=\sum_{\ell\not=j}u_{t\ell}^2+\e
\eea
and obtain
\bea
D&=&
\sum u_{1jt}^2 \sum_{k\not=j} {1+u_{kk}\over 1+u_{jj}}{u_{tk}^2\over u_{tt}(n+\Delta u)}
+
\sum_{j\not=k}{1\over u_{tt}(n+\Delta u)} u_{tk}u_{tk1}u_{tj}u_{tj1}
\nonumber\\
&&
+
\sum u_{1jt}^2 \sum_{k\not=j} {1+u_{kk}\over 1+u_{jj}}{1\over u_{tt}(n+\Delta u)}
(\sum_{\ell\not= j,k}u_{t\ell}^2+\e).
\nonumber\\
\eea
Thus we arrive at
\bea
D&=& {1\over 2}\sum_{j\not=k}
{1\over (1+u_{jj})(1+u_{kk})u_{tt}(n+\Delta u)}(u_{tj1}u_{tk}(1+u_{kk})+u_{tk1}u_{tj}(1+u_{jj}))^2
\nonumber\\
&&
+
\sum u_{1jt}^2 \sum_{k\not=j} {1+u_{kk}\over 1+u_{jj}}{1\over u_{tt}(n+\Delta u)}
(\sum_{\ell\not= j,k}u_{t\ell}^2+\e).
\eea

\subsection{A third formula for $Q$}

We summarize the expression for $Q$ obtained in this manner

\bea
Q&=&
{1\over 2}\sum_{j\not=k} \big\{({u_{tt}\over (1+u_{jj})(1+u_{kk})(n+\Delta u)})^{1\over 2}
((1+u_{kk})u_{jj1}-(1+u_{jj})u_{kk1})\nonumber\\
&&
+
({u_{tk}u_{tk1}\over 1+u_{kk}}-{u_{tj}u_{tj1}\over 1+u_{jj}})
({(1+u_{jj})(1+u_{kk})\over u_{tt}(n+\Delta u)})^{1\over 2}
\big\}^2
\nonumber\\
&&
+{1\over 2}\sum_{j\not=k}
{1\over (1+u_{jj})(1+u_{kk})u_{tt}(n+\Delta u)}
(u_{tj1}u_{tk}(1+u_{kk})+u_{tk1}u_{tj}(1+u_{jj}))^2
\nonumber\\
&&
+
\sum {1\over u_{tt}(n+\Delta u)}
u_{1jt}^2\sum_{k\not=j} {1+u_{kk}\over 1+u_{jj}}(
\sum_{\ell\not= j,k}u_{t\ell}^2+\e)
+E
\eea
where the term $E$ is defined by
\bea
E\equiv u_{tt}\sum_{j\not=k}
{u_{1jk}^2\over 1+u_{jj}}
-2
\sum_{j\not=k} u_{tk}u_{tj1}u_{1jk}{1\over 1+u_{jj}}.
\eea

\subsection{The $E$ terms}

We rewrite the $E$ term as follows,
\bea
E&=& {1\over 2}\sum_{j\not=k}u_{tt}({u_{1jk}^2\over 1+u_{jj}}
+
{u_{1kj}^2\over 1+u_{kk}})
-
\sum_{j\not=k}u_{1jk}({u_{tk}u_{tj1}\over 1+u_{jj}}+{u_{tj}u_{tk1}\over 1+u_{kk}})
\nonumber\\
&=&
{1\over 2}
\sum_{j\not=k}
{1\over (1+u_{jj})(1+u_{kk})}
\big\{u_{1jk}^2((1+u_{jj})+(1+u_{kk}))u_{tt}\nonumber\\
&&
-
2 u_{1jk}(u_{tk}u_{tj1}(1+u_{kk})+u_{tj}u_{tk1}(1+u_{jj}))\big\}.
\eea

\subsection{Differences when $\# G\geq 3$}

It is here that there seems to be a significant difference between the cases of
$\# G=2$ and $\#G>2$.
When $\#G=2$ and $\mu_0=0$, we have, modulo $O(\vp)$,
\bea
\label{Difficulty1}
(1+u_{jj})+(1+u_{kk})=n+\Delta u, \forall j\neq k \in G,
\eea
but not otherwise.

\subsection{Case $\# G=2$ }

When $\#G=2$, the argument can be completed as follows.
Using the fact
that $(1+u_{jj})+(1+u_{kk})=n+\Delta u$, we can write
\bea
E= {1\over 2}\sum_{j\not=k}{u_{1jk}^2(n+\Delta u)u_{tt}
-
2 u_{1jk}(u_{tk}u_{tj1}(1+u_{kk})+u_{tj}u_{tk1}(1+u_{jj}))\over(1+u_{jj})(1+u_{kk})}
\eea
and hence
\bea
E&=&
{1\over 2}
\sum_{j\not=k}{1\over (1+u_{jj})(1+u_{kk})}
\big\{u_{1jk}\sqrt{(n+\Delta u)u_{tt}}
-
{u_{tk}u_{tj1}(1+u_{kk})+u_{tj}u_{tk1}(1+u_{jj})
\over
\sqrt{u_{tt}(n+\Delta u)}}\big\}^2
\nonumber\\
&&
-
{1\over 2}\sum_{j\not=k}{1\over (1+u_{jj})(1+u_{kk})}
{1\over u_{tt}(n+\Delta u)}
(u_{tk}u_{tj1}(1+u_{kk})+u_{tj}u_{tk1}(1+u_{jj}))^2
\nonumber
\eea
Note now that the next to last term in $E$ cancels a term in the third formula
for $Q$.

\subsection{A fourth formula for $Q$}

Thus we have obtained, when $\#G=2$,
\bea
\label{Q-D1e}
Q&=&
{1\over 2}\sum_{j\not=k} \big\{({u_{tt}\over (1+u_{jj})(1+u_{kk})(n+\Delta u)})^{1\over 2}
((1+u_{kk})u_{jj1}-(1+u_{jj})u_{kk1})\nonumber\\
&&
+
({u_{tk}u_{tk1}\over (1+u_{kk})}-{u_{tj}u_{tj1}\over (1+u_{jj})})
({(1+u_{jj})(1+u_{kk})\over u_{tt}(n+\Delta u)})^{1\over 2}
\big\}^2
\nonumber\\
&&
+
{1\over 2}
\sum_{j\not=k}{\big\{u_{1jk}\sqrt{(n+\Delta u)u_{tt}}
-
{u_{tk}u_{tj1}(1+u_{kk})+u_{tj}u_{tk1}(1+u_{jj})
\over
\sqrt{u_{tt}(n+\Delta u)}}\big\}^2\over (1+u_{jj})(1+u_{kk})}
\nonumber\\
&&
+
\sum {1\over u_{tt}(n+\Delta u)}
u_{1jt}^2\sum_{k\not=j} {1+u_{kk}\over (1+u_{jj})}(
\sum_{\ell\not= j,k}u_{t\ell}^2+\e).
\eea
This expression for $Q$ shows that it is non-negative.
The proof of Proposition \ref{Dresult} is complete.

\section{Proof of the main theorems}
\setcounter{equation}{0}

Theorems \ref{1} is a consequence
of Proposition \ref{MAresult}.
In this case, since $\e>0$, we have $D_{tx}^2u+I_{n+1}>0$.
If the matrix $D_x^2u+I_n-\mu_0I_n$ has a zero
eigenvalue on the boundary of $X\times I$, there is nothing to prove.
Otherwise, if $x_0$ is an interior point with a zero eigenvalue of
multiplicity $K\geq 1$, then all possible values of $K$ are covered
by Proposition \ref{MAresult} when $n\leq 2$. Thus $D_x^2u+I_n-\mu_oI_n$
vanishes everywhere, and in particular again on the boundary.

\smallskip

The argument for Theorem \ref{2} is similar using Proposition \ref{Dresult},
except that we need to show first that
the solution $u$ satisfies the
space convexity condition $D^2_xu+I_n>0$ for each $t$.
In view of Proposition
\ref{Dresult}, we need to create a homotopic deformation path. Note that $u=1+t^2$ is the solution to
the equation
\[u_{tt}(n+\Delta u)-\sum_{k=1}^n u_{tk}^2=2n,\]
with boundary data
\[u(x,0)=1,
\quad
u(x,1)=2,\quad
 \forall x\in X.\]
For $\e>0$ and given boundary data $u^0, u^1$ satisfying (\ref{strictconvexity}),  for each $0\le s\le 1$, we consider the following family of equations
\begin{equation}
\label{Es}
u_{tt}(n+\Delta u)-\sum_{k=1}^n u_{tk}^2=s\e+2n(1-s),
\end{equation}
with boundary data
\begin{equation}
\label{Bs}
u(x,0)=su^0(x)+1-s,
\qquad
u(x,1)=su^1(x)+2(1-s),\quad \forall x\in X.
\end{equation}
It is easy to see that boundary data (\ref{Bs}) satisfies the condition (\ref{strictconvexity}) (possibly with different
$\lambda>0$, but independent of $s$). By \cite{H}, the
equation (\ref{Es}) has a
unique smooth solution for each $s\in [0,1]$. By continuity, the solution $u$ satisfies $D^2_xu+I_n>0$ for each $t$ when $s>0$ is small.
Let $s_0>0$ be the first value of $s$ where $D_x^2u+I_n$
has a zero eigenvalue for some $(x,t)\in X\times T$, if such a point
exists. Then for all $s<s_0$, $D_x^2u+I_n>0$ and we can apply
Proposition \ref{Dresult} to the equation with Dirichlet data corresponding to
this value of $s$. It follows that $D_x^2u+I_n\geq \lambda$ for $s<s_0$.
By continuity, this inequality still holds at $s=s_0$. This is a contradiction,
and thus no point with $D_x^2u+I_n$ with a zero eigenvalue exists.
This establishes the strict space convexity of the solution for all $0\leq s\leq 1$.
By applying again
Proposition \ref{Dresult}, we obtain the precise lower bound
$D^2_xu+I_n\ge \lambda$ for all $t, s\in [0,1]$.
Theorem \ref{2} is proved.

\smallskip
Finally, Theorem \ref{th2} follows directly from the last statement of Proposition \ref{Dresult}, and
Theorem \ref{th3} follows directly from
Theorem \ref{th2}.

\bibliographystyle{amsplain}

\begin{thebibliography}{99}

\bibitem{ALL} Alvarez, O., J.M. Lasry, and P.L. Lions,
``Convexity viscosity solutions and state constraints",
J. Maths. Pures Appl. 76 (1997) 265-288

\bibitem{BG1} Bian, B. and P. Guan,
``A microscopic convexity principle for nonlinear partial differential
equations", Inventiones Math. 177 (2009) 307-335.

\bibitem{BG2} Bian, B. and P. Guan,
``A structural condition for microscopic convexity principle",
Discrete and continuous dynamical ststems, 28, (2010) 789-807.

\bibitem{B} Blocki, Z., ``On geodesics in the space of K\"ahler metrics",
preprint, 2009.

\bibitem{BL} Brascamp, H.J. and E.H. Lieb,
``On extensions of the Brunn-Minkowski and Prokopa-Leindler
theorems, including inequalities for log-concave functions
with an application to the diffusion equation",
J. Funct. Anal. 22 (1976) 366-389.

\bibitem{CF}
Caffarelli, L. and A. Friedman ,
``Convexity of solutions of some semilinear elliptic equations",
Duke Math. J. 52 (1985) 431-455.

\bibitem{CGM} Caffarelli, L., P. Guan, and X. Ma,
``A constant rank theorem for solutions of fully nonlinear elliptic
equations", Comm. Pure Appl. Math. 60 (2007) 1769-1791.

\bibitem{C} Chen, X.X.,
``The space of K\"ahler metrics", J. Differential Geom.
56 (2000) 189-234.

\bibitem{CH} Chen, X.X. and W. He,
``The space of volume forms", arXiv: 0810.3880

\bibitem{CS} Chen, X.X. and S. Sun,
``Space of K\"ahler metrics (V): K\"ahler quantization",
arXiv: 0902.4149

\bibitem{CNS} Caffarelli, L., L. Nirenberg, and J. Spruck,
``The Dirichlet problem for non-linear second-order elliptic equations I.
Monge-Amp\`ere equations", Comm. Pure Appl. Math. XXXVII
(1984) 369-402.

\bibitem{D1} Donaldson, S.K.,
``Symmetric spaces, K\"ahler geometry, and Hamiltonian dynamics",
Amer. Math. Soc. Transl. 196 (1999) 13-33.

\bibitem{D2} Donaldson, S.K.,
``Nahm's equations and free boundary problems",
arXiv: 0709.0184

\bibitem{Gb} Guan, B.,
``The Dirichlet problem for Hessian equations on Riemannian manifolds".
Calc. Var. PartialDifferential Equations, 8 (1999), 45–69.

\bibitem{Gd} Guan, D.,
``On modified Mabuchi functional and Mabuchi moduli space
of K\"ahler metrics on toric bundles", Math. Res. Lett. 6 (1999)
no. 5-6, 547-555.


\bibitem{GLZ} Guan, P., Li, Q. and X. Zhang, `` A uniqueness theorem in K\"ahler geometry",
Math. Ann. Vol. 345, (2009) 377-393.

\bibitem{H} He, W.,
``The Donaldson equation", arXiv: 0810.4123

\bibitem{KL} Korevaar, N.J. and J. Lewis,
``Convex solutions of certain elliptic equations have constant rank Hessians",
Arch. Rat. Mech. Anal. 97 (1987) 19-32.

\bibitem{Lm} Lempert, L., ``La m\'etrique de Kobayashi et la repr\'esentation des domaines sur la boule". Bull. Soc. Math. France 109 (1981), 427-474.

\bibitem{L} Li, Q., ``Constant rank theorem in complex variables". Indiana Univ. Math. J. 58 (2009), 1235--1256.

\bibitem{M} Mabuchi, T.,
``Some symplectic geometry on compact K\"ahler manifolds",
Osaka J. Math. 24 (1987) 227-252.

\bibitem{PS1} Phong, D.H. and J. Sturm,
``The Monge-Amp\`ere operator and geodesics in the space of
K\"ahler potentials", Inventiones Math. 166 (2006) 125-149.

\bibitem{PS2} Phong, D.H. and J. Sturm,
``The Dirichlet problem for degenerate complex Monge-Amp\`ere equations",
arXiv: 0904.1898, to appear in Commun. Anal. Geometry.

\bibitem{PS3} Phong, D.H. and J. Sturm,
``Test configurations for K-stability and geodesic rays", J. Symplectic Geom. {\bf 5} (2007) no. 2,
221-247.

\bibitem{S} Semmes, S.,
``Complex Monge-Amp\`ere equations and symplectic manifolds",
Amer. J. Math. 114 (1992) 495-550.

\bibitem{Y} Singer, I.M., B. Wong, S.S.T. Yau, and S.T. Yau,
``An estimate of gap of the first two eigenvalues of the Schr\"odinger
operator", Ann. Scuola Norm. Sup. Pisa Cl. Sci. (4) 12 (1985) 319-333.

\end{thebibliography}

\end{document}